\documentclass[12pt]{amsart}
\usepackage{amsfonts}
\usepackage{amsmath}
\usepackage{amssymb}
\usepackage{latexsym}

\newcommand{\bR}{\mathbb{R}}

\newtheorem{theorem}{Theorem}
\theoremstyle{plain}
\newtheorem{corollary}[theorem]{Corollary}
\newtheorem{definition}[theorem]{Definition}
\newtheorem{proposition}[theorem]{Proposition}

\newtheorem{remark}[theorem]{Remark}
\numberwithin{equation}{section}
\renewcommand\endproof{\hfill $\Box$\vskip 0.15in}

\begin{document}
\title[Second-Order SPDEs]{Stochastic Parabolic Equations of Full Second
Order}
\author{S. V. Lototsky}
\curraddr[S. V. Lototsky]{Department of Mathematics, USC\\
Los Angeles, CA 90089}
\email[S. V. Lototsky]{lototsky@math.usc.edu}
\urladdr{http://www-rcf.usc.edu/$\sim$lototsky}
\author{B. L. Rozovskii}
\curraddr[B. L. Rozovskii]{Division of Applied Mathematics \\
Brown University\\
Providence, RI 02912}
\email[B. L. Rozovskii]{rozovsky@dam.brown.edu}
\thanks{S. V. Lototsky acknowledges support from NSF
 CAREER award
DMS-0237724. B. L. Rozovskii acknowledges support
 from NSF Grant DMS
0604863, ARO Grant W911NF-07-1-0044, and
ONR Grant N00014-07-1-0044}

\begin{abstract}
A procedure is described for defining a generalized solution for
stochastic differential equations using the Cameron-Martin version
of the Wiener Chaos expansion. Existence and uniqueness of this
Wiener Chaos solution is established for parabolic stochastic PDEs
such that both the drift and the diffusion operators are of the
second order.
\end{abstract}

\maketitle

\section{Introduction}

Consider a stochastic evolution equation
\begin{equation}
du(t)=({\mathcal{A}}u(t)+f(t))dt+({\mathcal{M}}u(t)+g(t))dW(t),
\label{eq:intr}
\end{equation}%
where ${\mathcal{A}}$ and ${\mathcal{M}}$ are differential operators, and $W$
is a Wiener process on a probability
space $(\Omega, \mathcal{F},\mathbb{P})$.
 Traditionally, this equation was studied under the
following assumptions:
\begin{itemize}
\item \textbf{A1. }The operator ${\mathcal{A}}$ is elliptic, the order of
the operator ${\mathcal{M}}$ is less than the order of ${\mathcal{A}}$, and
${\mathcal{A}}-\frac{1}{2}{\mathcal{MM}}^{\star }$ is elliptic (possibly
degenerate) operator,
\end{itemize}
In fact, it is well known that unless assumption \textbf{A1} holds, equation
(\ref{eq:intr}) has no solutions in $L_{2}(\Omega ;X)$ for any reasonable
choice of the state space $X.$

It was shown recently (see \cite{LR_shir, LR_AP, MR_WCSIAM}
 and the references therein) that if only
the operator ${\mathcal{A}}$ is elliptic and the order of
${\mathcal{M}}$ is \emph{smaller} than the order of $A,$\ then
there exists a unique generalized solution of equation
(\ref{eq:intr}). This solution is often referred to as
Wiener Chaos solution. $\ $It is given by the Wiener chaos expansion $%
u\left( t\right) =\sum_{\left\vert \alpha \right\vert <\infty }u_{\alpha
}\left( t\right) \xi _{\alpha }$, where $\left\{ \xi _{\alpha }\right\}
_{\left\vert \alpha \right\vert <\infty }$ is the Cameron-Martin orthonormal
basis in the space $L_{2}(\Omega ;\mathcal{F}^{W};X)$ of square
integrable random elements in $X$ measurable with respect to
 to the sigma-algebra $\mathcal{F}^{W}$ generated by the
Wiener process.
 The Cameron-Martin basis $\left\{ \xi _{\alpha }\right\} $ is
indexed by multiindices $\alpha =\left( \alpha _{1},\alpha _{2},...\right) .$
It was shown that for certain positive weights $Q=\left\{ q\left( \alpha
\right) \right\} _{\left\vert \alpha \right\vert <\infty }$, the weighted
norm
\begin{equation*}
\left\Vert u\right\Vert _{Q,X}^{2}:=\sum_{\left\vert \alpha \right\vert
<\infty }q^{2}\left( \alpha \right) \left\Vert u_{\alpha }\right\Vert
_{L_{2}\left( (0,T);X\right) }^{2}<\infty ,
\end{equation*}%
where $X$ is the appropriate Hilbert space characterizing the
\textquotedblleft regularity\textquotedblright\ of
the solution. Note that without assumption {\bf A1}
\begin{equation*}
E\left\Vert u\right\Vert _{L_{2}\left( (0,T);X\right) }^{2}=\sum_{\left\vert
\alpha \right\vert <\infty }\left\Vert u_{\alpha }\left( t\right)
\right\Vert _{L_{2}\left( (0,T);X\right) }^{2}=\infty .
\end{equation*}

In this paper, we consider the Cauchy problem for the following stochastic
partial differential equation:
\begin{equation}
du=(a_{ij}D_{i}D_{j}u+b_{i}D_{i}u+cu+f)dt+(\rho _{ij}D_{i}D_{j}u+\sigma
_{i}D_{i}u+\nu u+g)dW,\ t\in \lbrack 0,T],\ x\in \mathbb{R}.
\label{eq: spdemain}
\end{equation}%
In contrast to the previous work, this is a parabolic SPDE of the \emph{full
}second order, in that the drift and diffusion operators have the same order
$2$. We construct a scale of weighted Wiener chaos spaces (related
but not identical to Kondratiev's spaces) and prove the existence and
uniqueness of the solution in the spaces from this scale.

\section{Constructing a Solution: an Example}

Let ${\mathbb{F}} =(\Omega, {\mathcal{F}}, \{{\mathcal{F}}_{t}\}_{0\leq
t\leq T}, {\mathbb{P}})$ be a stochastic basis with the usual assumptions
and $W=W(t),\;0\leq t\leq T$, a standard Wiener process on ${\mathbb{F}}.$
For a Hilbert space $X$, denote by $L_2(\mathbb{W};X)$ the collection of $X$%
-valued random elements that are square integrable ($\mathbb{E}%
\|\cdot\|_{X}^2 <\infty$) and are measurable with respect to the
sigma-algebra generated by $W(t), t\in [0,T]$.

Consider the It\^{o}\ equation
\begin{equation}  \label{eq:ex1}
u(t,x)=e^{-x^2/2}+\int_0^tu_{xx}(s,x)ds+ \int_0^t u_{xx}(s,x)dW(s),\ t\in
[0,T],\ x\in \mathbb{R}.
\end{equation}
If there is a solution, its Fourier transform in space, $\widehat{u}(t,y)=(1/%
\sqrt{2\pi})\int_{\mathbb{R}} e^{-ixy}u(t,x)dx$ satisfies
\begin{equation}  \label{eq:ex2}
\widehat{u}(t,y)=e^{-y^2/2} -y^2\int_0^t\widehat{u}(s,y)ds -y^2\int_0^t
\widehat{u}(s,y)dW(t), \ t\in [0,T],\ y\in \mathbb{R}.
\end{equation}
For each fixed $y$, \eqref{eq:ex2} defines a geometric Brownian motion:
\begin{equation}  \label{eq:ex3}
\widehat{u}(t,y)=e^{-(1+t)y^2-(y^4/2)t-y^2W(t)}.
\end{equation}
Let $H^{\gamma}(\mathbb{R})$ be the Sobolev space
\begin{equation}
\left\{ f: \int_{\mathbb{R}}
(1+|y|^2)^{\gamma} |\widehat{f}(y)|^2dy<
\infty\right\}.
\end{equation}
Since
\begin{equation}  \label{eq:ex4}
\mathbb{E}|\widehat{u}(t,y)|^2 =e^{-2(1+t)y^2+y^4t},
\end{equation}
the solution of \eqref{eq:ex1} cannot be an element of $L_2\big(\mathbb{W};
L_2((0,T);H^{\gamma}(\mathbb{R}))\big)$ for \emph{any} $\gamma\in \mathbb{R}$%
, even though the initial condition is non-random and is an element of $%
H^{\gamma}(\mathbb{R})$ for \emph{every} $\gamma\in \mathbb{R}$.

Let us try another approach. Once again, assuming that the solution exists,
we apply the It\^{o}\ formula to the product $u(t,x)\mathcal{E}_h(t)$, where
\begin{equation}  \label{eq:cE}
\mathcal{E}_h(t)=\exp\left(\int_0^th(s)dW(s) -\frac{1}{2}\int_0^th^2(s)ds
\right)
\end{equation}
and $h=h(t)$ is a smooth deterministic function. Since
\begin{equation}  \label{eq:cE1}
\mathcal{E}_h(t)=1+\int_0^t\mathcal{E}_h(s)h(s)dW(s),
\end{equation}
we conclude that the function
\begin{equation}  \label{eq:Str}
u_h(t,x)=\mathbb{E}\big(u(t,x)\mathcal{E}_h(t)\big),
\end{equation}
if defined, must satisfy the heat equation
\begin{equation}  \label{eq:ex10}
u_h(t,x)=e^{-x^2/2}+\int_0^t(1+h(s)) \frac{\partial^2 u_h(s,x)}{\partial x^2}%
ds.
\end{equation}
If $\sup_t|h(t)|<1$, then this equation has a unique solution in every $%
H^{\gamma}(\mathbb{R})$ and
\begin{equation}  \label{eq:ex11}
u_h(t,x)=\mathbb{E}\exp(-(X(t,x))^2/2),
\end{equation}
where
\begin{equation}  \label{eq:ex12}
X(t,x)=x+\int_0^t\sqrt{2(1+h(s))}\,dW(s).
\end{equation}
In other words, while existence of a solution of equation \eqref{eq:ex1} is
still unclear, we now have a family of functions $u_h(t,x)$ defined by %
\eqref{eq:ex11}. All we need now is a systematic procedure of relating the
family of deterministic functions $u_h=u_h(t,x)$ to a random process $%
u=u(t,x)$; then this process is natural to call a solution of \eqref{eq:ex1}.

Here is a possible way of constructing a stochastic process from $u_h$. Let $%
\mathfrak{m}= \{m_k,\ k\geq 1\}$ be the Fourier cosine basis in $L_2((0,T))$%
:
\begin{equation}  \label{basis}
m_1(s)\!=\!\frac 1{\sqrt{T}};\ m_k(t)\!= \!\sqrt{\frac{2}{T} } \cos \left(
\frac{\pi (k-1) t}{T} \right), \, k>1; \ 0\leq t \leq T.
\end{equation}
Then
\begin{equation}  \label{ex:eq21}
h(t)=\sum_{k\geq 1}h_km_k(t),
\end{equation}
For every fixed $t\in [0,T]$ and $\gamma\in \mathbb{R}$, we can now
interpret the function $u_h(t,\cdot)$ as a mapping from the set of sequences
$h=(h_1, h_2,\ldots)$ to the space $H^{\gamma}(\mathbb{R}^d)$, and, as
equalities \eqref{eq:ex11} and \eqref{eq:ex12} suggest, this mapping is
analytic in the region $\{h: \sum_{k\geq 1}h_k^2<\varepsilon\}$ for
sufficiently small $\varepsilon$. We will now compute the derivatives of
this mapping.

Let $\mathcal{J}$ be the collection of multi-indices $\alpha=\{\alpha_{k},\
k\geq 1\}$. Each $\alpha\in \mathcal{J}$ has non-negative integer elements $%
\alpha_{k}$ and
\begin{equation}  \label{eq:multin1}
|\alpha|=\sum_{k}\alpha_{k}<\infty.
\end{equation}
We also use the notation
\begin{equation}  \label{eq:multin2}
\alpha!=\prod_{k}\alpha_{k}!
\end{equation}
and consider special multi-indices, $\alpha=(0)$ with $|\alpha|=0$ and $%
\alpha=\varepsilon_i$, with $|\alpha|=1$, $\alpha_i=1$.

For each $\alpha \in \mathcal{J}$ define
\begin{equation}  \label{eq:ex22}
u_{\alpha}(t,x)= \frac{1}{\sqrt{\alpha!}} \frac{\partial^{|\alpha|}u_{h}(t,x)%
} {\partial h_1^{\alpha_1}\partial h_2^{\alpha_2} \cdots} \Bigg|_{h=0}.
\end{equation}
Then
\begin{equation}  \label{eq:ex23}
u_h(t,x)=\sum_{\alpha\in \mathcal{J}} u_{\alpha}(t,x) \frac{h^{\alpha}}{%
\sqrt{\alpha!}},
\end{equation}
where
\begin{equation}
h^{\alpha}=\prod_{k\geq 1}h_k^{\alpha_k}.
\end{equation}
On the other hand, by direct computation,
\begin{equation}
\mathcal{E}_h(t)=\mathbb{E}(\mathcal{E}_h(T)|\mathcal{F}^W_t)=\sum_{\alpha%
\in \mathcal{J}} \frac{h^{\alpha}}{\sqrt{\alpha!}}\xi_{\alpha}(t),
\end{equation}
where
\begin{equation}  \label{eq:xial}
\xi_{\alpha}(t)=\mathbb{E}(\xi_{\alpha}|\mathcal{F}^W_t),\ \xi_{\alpha}=%
\frac{1}{\sqrt{\alpha!}} \prod_{k\geq 1} H_{\alpha_k}\left(\int_0^T
m_k(t)dW(t)\right),
\end{equation}
and
\begin{equation}
H_n(x)=(-1)^ne^{x^{2}/2} \frac{d^{n}}{dx^{n}}e^{-x^{2}/2}
\end{equation}
is $n$-th Hermite polynomial. It is a standard fact \cite{CM} that the
collection $\{\xi_{\alpha}, \ \alpha \in \mathcal{J}\}$ is an orthonormal
basis in $L_2(\mathbb{W};\mathbb{R})$.

The functions $u_{\alpha}(t,x),\ \alpha \in \mathcal{J},$ uniquely determine
$u_h(t,x)$ according to (\ref{eq:ex23}). On the other hand, if
\begin{equation}  \label{eq:conv}
\sum_{\alpha\in \mathcal{J}} \|u_{\alpha}(t)\|_{H^{\gamma}(\mathbb{R}%
)}^2<\infty,
\end{equation}
then the $H^{\gamma}(\mathbb{R})$-valued random process
\begin{equation}  \label{eq:proc}
u(t,x)=\sum_{\alpha\in \mathcal{J}} u_{\alpha}(t,x)\xi_{\alpha}
\end{equation}
satisfies $\mathbb{E}(u(t,x)\mathcal{E}_h(t))=u_h(t,x)$; if, in addition, $u$
is $\mathcal{F}^W_t$-adapted, then also
\begin{equation}  \label{eq:formal}
u(t,x)=\sum_{\alpha\in \mathcal{J}} u_{\alpha}(t,x)\xi_{\alpha}(t).
\end{equation}
If condition \eqref{eq:conv} fails, then \eqref{eq:proc} is a formal series,
which we \emph{define} to be the stochastic process corresponding to the
family $u_h$.

As \eqref{eq:ex4} suggests, if $u_h$ is the solution of \eqref{eq:ex10},
then \eqref{eq:conv} fails for every $\gamma$. Let us now see how fast the
series diverges. Equality \eqref{eq:ex10} implies
\begin{equation}  \label{s-system-ex}
\begin{split}
u_{(0)}(t,x)&=e^{-x^2/2}+\int_0^t \frac{\partial^2u_{(0)}(s,x)}{\partial x^2}%
ds, \ |\alpha|=0; \\
u_{\epsilon_{i}}(t)&= \int_0^t \frac{\partial^2u_{\epsilon_i}(s,x)}{\partial
x^2}ds+ \int_0^t \frac{\partial^2u_{(0)}(s,x)}{\partial x^2}m_i(s)ds, \
|\alpha|=1; \\
u_{\alpha}(t)&= \int_0^t \frac{\partial^2u_{\alpha}(s,x)}{\partial x^2}ds +
\sum_{k=1}^{\infty}\sqrt{\alpha_{k}} \int_0^t\frac{\partial^2
u_{\alpha-\epsilon_{k}}(s,x)} {\partial x^2} m_k(s)ds, \ |\alpha|>1.
\end{split}%
\end{equation}
Equations of the type \eqref{s-system-ex} have been studied \cite[Section 6
and References]{LR_shir}. In particular, it is known that
\begin{equation}  \label{eq:ex231}
\sum_{|\alpha|=n} \|u_{\alpha}(t)\|_{H^{\gamma}(\mathbb{R})}^2 = \frac{t^n}{%
n!} \|D^{2n}_x\Phi_t u_0\|_{H^{\gamma}(\mathbb{R})}^2,
\end{equation}
where $D_x=\partial/\partial x$, $\Phi_t$ is the heat semigroup, and $%
u_0(x)=e^{-x^2/2}$. To simplify further computation, let us assume that $%
\gamma=0$. Then, switching to the Fourier transform,
\begin{equation}
\|D^{2n}_x\Phi_t u_0\|_{L_2(\mathbb{R}^d)}^2= \int_{\mathbb{R}%
}|y|^{4n}e^{-y^2(t+1)}dy=
\frac{\Gamma\left(2n+\frac{1}{2}\right)}{(1+t)^{2n}.%
}
\end{equation}
Using Stirling's formula for the Gamma function $\Gamma$,
\begin{equation}  \label{eq:ex24}
\sum_{|\alpha|=n} \|u_{\alpha}(t)\|_{L_2(\mathbb{R})}^2 =\left(\frac{2\sqrt{t%
}}{1+t}\right)^{2n} C(n) n!,
\end{equation}
where the numbers $C(n)$ are uniformly bounded from above and below. Similar
result holds in every $H^{\gamma}(\mathbb{R})$. Thus, \eqref{eq:conv} does
not hold, but instead, by (\ref{eq:ex24}), we have
\begin{equation}
 \label{eq:ex25}
\sum_{\alpha\in \mathcal{J}} \frac{1}{|\alpha|!}
 \|u_{\alpha}(t)\|_{H^{\gamma}(%
\mathbb{R})}^2 <\infty.
\end{equation}
 We denote by $(\mathfrak{L})_{0,0}(%
\mathbb{W};H^{\gamma}(\mathbb{R}))$ the collection of formal series %
\eqref{eq:formal} satisfying \eqref{eq:ex25};
 the reason for using $(\mathfrak{L})_{0,0}$ in the
notation will become clear later. Note that we had equalities in all
computations for equation \eqref{eq:ex1} that lead to \eqref{eq:ex25}, which
suggests that $(\mathfrak{L})_{0,0}(\mathbb{W};H^{\gamma}(\mathbb{R}))$ is
the natural solution space for equation \eqref{eq:ex1}. For a more general
stochastic parabolic equation of full second order in $\mathbb{R}^d $, the
natural solution space turns out to be $(\mathfrak{L})_{p,q}(\mathbb{W}%
;L_2((0,T);H^{\gamma}(\mathbb{R}^d)))$ for suitable $p,q\leq 0$.

In the next section we address the following questions:

\begin{enumerate}
\item How to define the spaces $(\mathfrak{L})_{p,q}(\mathbb{W}; X)$ for $%
p,q\in \mathbb{R}$ without relying on an orthonormal basis in $L_2((0,T))$?

\item How to construct a solution of a general stochastic parabolic
equations of full second order?
\end{enumerate}

\section{General Constructions and the Main Result}

As before, let ${\mathbb{F}} =(\Omega, {\mathcal{F}}, \{{\mathcal{F}}%
_{t}\}_{0\leq t\leq T}, {\mathbb{P}})$ be a stochastic basis with the usual
assumptions and $W=W(t),\;0\leq t\leq T$, a standard Wiener process on ${%
\mathbb{F}}.$ Denote by $H^s=H^s((0,T))$, $s\geq 0$, the Sobolev spaces on $%
(0,T)$ with norm
\begin{equation}
\|\cdot\|_s=\|\Lambda^{s/2}\cdot\|_{L_2((0,T))},
\end{equation}
where $\Lambda$ is the operator
\begin{equation}  \label{eq:operT}
1-\frac{T^2}{\pi^2}\,\frac{d^2}{dt^2}
\end{equation}
with Neumann boundary conditions. This norm extends to functions of several
variables via the tensor product of the spaces $H^s$.

\begin{definition}
Given real numbers $p,q$ and a Hilbert space $X$, $(\mathfrak{L})_{p,q}(%
\mathbb{W};X)$ is the closure of the set of $X$-valued random elements
\begin{equation}  \label{eq:spS}
\eta=\eta_0+\sum_{k=1}^N \int_0^T\int_0^{s_{k}}\ldots\int_0^{s_2}
\eta_k(s_1,\ldots,s_k)dW(s_1)\ldots dW(s_{k-1})dW(s_k), \ N\geq 1,
\end{equation}
with respect to the norm
\begin{equation}  \label{eq:spS1}
\|\eta\|_{p,q;X}^2= \|\eta_0\|_X^2+\sum_{k=1}^N \frac{2^{kp}}{(k!)^2}\,\| \,
\|\eta_k\|_{q}\,\|_{X}^2,
\end{equation}
where each $\eta_k$, $k\geq 1$, is a smooth symmetric function from $[0,T]^k$
to $X$.
\end{definition}

\begin{remark}
(a) It is known \cite{HOUZ, Nualart} that, for $\eta$ of the type %
\eqref{eq:spS},
\begin{equation}  \label{eq:L2norm}
\mathbb{E}\|\eta\|_X^2=\|\eta_0\|_X^2+ \sum_{k=1}^N \frac{1}{k!}\| \,
\|\eta_k\|_{0}\,\|_{X}^2.
\end{equation}

(b) The definition of each individual $(\mathfrak{L})_{p,q}(\mathbb{W};X)$
inevitably involves arbitrary choices, such as the norm
in $H^q((0,T))$. Further analysis shows that
different choices result in shifts of the indices $p,q$, and
the space $\cup_{p,q}(\mathfrak{L})_{p,q}(\mathbb{W};X)$ does not depend on
any arbitrary choices. In the white noise setting, where $\Omega$ is the
space $\mathcal{S}'(\bR^d)$ of the Schwartz distributions and $\mathbb{P}$ is the
normalized Gaussian measure on $\mathcal{S}$, the inductive limit $%
\cup_{p,q}(\mathfrak{L})_{p,q}(\mathbb{W};\mathbb{R})$ is the Kondratiev
space $(\mathcal{S})_{-1}$ \cite{HOUZ}.
\end{remark}

\begin{remark}
If $X$ is the Sobolev space $H^{\gamma}(\mathbb{R}^d)$, then we denote the
norm $\|\cdot\|_{p,q;X}$ by $\|\cdot\|_{p,q;\gamma}$:
\begin{equation}
\|\cdot\|_{p,q;H^{\gamma}(\mathbb{R}^d)}= \|\cdot\|_{p,q;\gamma}.
\end{equation}
\end{remark}

\begin{proposition}
Let $\eta=f\xi_{\alpha}$, where $f\in X$ and $\xi_{\alpha}$ is defined by %
\eqref{eq:xial}. Then
\begin{equation}  \label{eq:norm-b}
\|\eta\|_{p,q;X}^2= \frac{2^{|\alpha|p}}{|\alpha|!} \mathbb{N}%
^{2q\alpha}\|f\|^2_X,
\end{equation}
where
\begin{equation}
\mathbb{N}^{2q\alpha}=\prod_{k\geq 1} k^{2q\alpha_k}.
\end{equation}
\end{proposition}

\textbf{Proof.} Let $|\alpha|=n$. It is known \cite{Ito} that
\begin{equation}
\xi_{\alpha}= \frac{1}{\sqrt{\alpha!}} \int_0^T\int_0^{s_n} \ldots
\int_0^{s_2} E_{\alpha}(s_1,\ldots,s_n) dW(s_1)\ldots dW(s_{n-1})dW(s_n),
\end{equation}
where $E_{\alpha}$ is the symmetric function
\begin{equation}  \label{eq:Ealpha}
E_{\alpha}(s_1,\ldots,s_n)= \sum_{\sigma\in \mathcal{P}_n }
m_{i_1}(s_{\sigma(1)})\ldots m_{i_{n}}(s_{\sigma(n)}).
\end{equation}
In (\ref{eq:Ealpha}), the summation is over all permutations of $%
\{1,\ldots,n\}$, the functions $m_k$ are defined in \eqref{basis}, and the
positive integer numbers $i_1\leq i_2\leq \ldots i_n$ are such that, for
every sequence $(b_k,k\geq 1)$ of positive numbers,
\begin{equation}
\prod_{k\geq 1}b_k^{\alpha_k}= b_{i_1}\cdot b_{i_2}\cdot \ldots \cdot
b_{i_n}.
\end{equation}
For example, if $\alpha=(1,0,2,0,0,4,0,0,\ldots)$, then $|\alpha|=7$ and $%
i_1=1$, $i_2=i_3=3$, $i_4=\ldots=i_7=6$. Thus, in the notations of %
\eqref{eq:spS1}, we have
\begin{equation}
\eta_k=
\begin{cases}
\displaystyle \frac{1}{\sqrt{\alpha!}}E_{\alpha}\, f, & \mathrm{if} \ k=n,
\\
0, & \mathrm{otherwise}%
\end{cases}%
\end{equation}
Note that
\begin{equation}
\|E_{\alpha}\|_0=\sqrt{\alpha!}\sqrt{n!}.
\end{equation}

By definition \eqref{eq:operT} of the operator $\Lambda$ we have
\begin{equation}  \label{eq:eigT}
\Lambda^{q/2} m_k=k^{q} m_k \ \ \ \mathrm{or}\ \ \ \|E_{\alpha}\|_{q}^2= {%
\mathbb{N}^{2q\alpha}\alpha!n!}.
\end{equation}
The result now follows. \endproof

\begin{corollary}
A formal series
\begin{equation}
\eta=\sum_{\alpha\in \mathcal{J}} \eta_{\alpha}\xi_{\alpha},
\end{equation}
with $\eta_{\alpha}\in X$, is an element of $(\mathfrak{L})_{p,q}(\mathbb{W}%
;X)$ if and only if
\begin{equation}  \label{eq:pqnorm-ser}
\sum_{\alpha\in \mathcal{J}} \frac{2^{p|\alpha|}\mathbb{N}^{2q\alpha}}{%
|\alpha|!} \|\eta_{\alpha}\|^2_X< \infty.
\end{equation}
\end{corollary}

\textbf{Proof.} This follows from \eqref{eq:eigT} and the equality
\begin{equation}
\|E_{\alpha}+E_{\beta}\|_0^2=\|E_{\alpha}\|_0^2+ \|E_{\beta}\|_0^2,\
\alpha\not=\beta.
\end{equation}
\endproof

Denote by $(\mathfrak{L})^{p,q}(\mathbb{W})$ the Hilbert space dual of $(%
\mathfrak{L})_{p,q}(\mathbb{W};\mathbb{R})$ relative to the inner product in
$L_2(\mathbb{W};\mathbb{R})$, and by $\langle\!\langle\cdot,\cdot\rangle\!%
\rangle$ the corresponding duality.
In the white noise setting, $\cap_{p,q}(%
\mathfrak{L})^{p,q}(\mathbb{W})$ is the space $(\mathcal{S})_{1}$ of the
Kondratiev test functions \cite{HOUZ}. If $\eta\in(\mathfrak{L})_{p,q}(%
\mathbb{W};X)$ and $\zeta\in (\mathfrak{L})^{p,q}(\mathbb{W})$, then $%
\langle\!\langle\eta,\zeta\rangle\!\rangle$ is defined and belongs to $X$.

For $h\in L_2((0,T))$, define
\begin{equation}  \label{eq:SE}
\mathcal{E}_h=\mathcal{E}_h(T)=\exp\left( \int_0^T h(s)dW(s)-\frac{1}{2}%
\int_0^T |h(s)|^2ds \right).
\end{equation}

\begin{proposition}
\label{prop:EhSpq} The random variable $\mathcal{E}_h$ is an element of $(%
\mathfrak{L})^{p,q}(\mathbb{W})$ if an only if
\begin{equation}  \label{eq:smallh}
\|h\|_{-q}^2< 2^{p}.
\end{equation}
\end{proposition}

\textbf{Proof.} Since
\begin{equation}
\mathcal{E}_h(T)=1+\int_0^Th(t)\mathcal{E}_h(t)dt,
\end{equation}
it follows that
\begin{equation}
\mathcal{E}_h=1+\sum_{k=1}^{\infty} \int_0^T\int_0^{s_k}\ldots \int_0^{s_2}
h(s_k)\cdots h(s_1) dW(s_1)\cdots dW(s_{k-1})dW(s_k).
\end{equation}
By \eqref{eq:spS1} and \eqref{eq:L2norm}, $\mathcal{E}_h\in (\mathfrak{L}%
)^{p,q}(\mathbb{W})$ if and only if
\begin{equation}
\sum_{k=1}^{\infty} \Big(2^{-p}\|h\|_{-q}^2\Big)^k < \infty,
\end{equation}
that is, $\|h\|_{-q}^2<2^p$. \endproof

\begin{definition}
We say that the function $h$ is sufficiently small if \eqref{eq:smallh}
holds for sufficiently large (positive) $-p,-q$.
\end{definition}

\begin{proposition}
If $u\in \bigcup_{p,q} (\mathfrak{L})_{p,q}(\mathbb{W};X)$ and $h$ is
sufficiently small, then
\begin{equation}
u_h=\langle\!\langle u,\mathcal{E}_h\rangle\!\rangle
\end{equation}
is an $X$-valued analytic function of $h$.
\end{proposition}

\textbf{Proof.} For every $u\in \bigcup_{p,q} (\mathfrak{L})_{p,q}(\mathbb{W}%
;X)$, there exist $p,q$ such that $u\in (\mathfrak{L})_{p,q}(\mathbb{W};X)$;
by Proposition \ref{prop:EhSpq}, $u_h$ will indeed be defined for
sufficiently small $h$. Similar to \eqref{eq:ex23} we have
\begin{equation}
u_h=\sum_{\alpha\in \mathcal{J}} \frac{u_{\alpha}h^{\alpha}} {\sqrt{\alpha!}}
\end{equation}
and this power series in $h^{\alpha}$ converges in some
(infinite-dimensional) neighborhood of zero. \endproof

From now on, $D_i=\partial /\partial x_i$, and the summation convention is
in force: $c_id_i=\sum_i c_id_i$, etc.

Consider the linear equation in ${\mathbb{R}}^{d}$
\begin{equation}  \label{eq:linG}
du=(a_{ij}D_{i}D_{j}u +b_{i}D_{i}u+cu+f)dt +(\rho_{ij}D_iD_ju
+\sigma_{i}D_{i}u +\nu u+g)dW
\end{equation}
with initial condition $u(0,x)=v(x)$, under the following \textbf{%
assumptions:}

\begin{enumerate}
\item[\textbf{{B0}}] All coefficients are non-random.

\item[\textbf{{B1}}] The functions $a_{ij}=a_{ij}(t,x)$, $%
\rho_{ij}=\rho_{ij}(t,x)$ are measurable and bounded in $(t,x)$ by a
positive number $C_0$, and

\begin{enumerate}
\item[(i)]
\begin{equation*}
|a_{ij}(t,x)-a_{ij}(t,y)| +|\rho_{ij}(t,x)-\rho_{ij}(t,y)| \leq C_0|x-y|,\
x,y\in{\mathbb{R}}^{d}, \ 0\leq t\leq T;
\end{equation*}

\item[(ii)] the matrix $(a_{ij})$ is uniformly positive definite, that is,
there exists a $\delta>0$ so that, for all vectors $y\in{\mathbb{R}}^{d}$
and all $(t,x)$, $a_{ij}y_{i}y_{j}\geq\delta|y|^{2}$.
\end{enumerate}

\item[\textbf{{B2}}] The functions $b_{i}=b_{i}(t,x)$, $c=c(t,x)$, $%
\sigma_{i}=\sigma_{i}(t,x)$, and $\nu=\nu(t,x)$ are measurable and bounded
in $(t,x)$ by the number $C_0$.

\item[\textbf{{B2}}]
\begin{equation}
u_0\in \bigcup_{p,q}(\mathfrak{L})_{p,q}(\mathbb{W};L_2(\mathbb{R}^d)),\
f,g\in \bigcup_{p,q} (\mathfrak{L})_{p,q}(\mathbb{W};L_2((0,T);H^{-1}(%
\mathbb{R}^d))).
\end{equation}
\end{enumerate}

For simplicity, we introduce the following notations for the differential
operators in \eqref{eq:linG}:
\begin{equation}
\mathcal{A}=a_{ij}D_{i}D_{j} +b_{i}D_{i}+c,\ \ \mathcal{B}=\rho_{ij}D_iD_j
+\sigma_{i}D_{i} +\nu.
\end{equation}

\begin{definition}

A solution $u$ of \eqref{eq:linG} is an element of $\bigcup_{p,q} (\mathfrak{%
L})_{p,q}(\mathbb{W};L_2((0,T);H^{1}(\mathbb{R}^d)))$ such that, for all
sufficiently small $h$ and all $t\in [0,T]$, the equality
\begin{equation}  \label{eq:Str-sol}
u_h(t,x)=v_{h}(x)+\int_0^t(\mathcal{A}+h(s)\mathcal{B}) u_h(s,x)ds
\end{equation}
holds in $H^{-1}(\mathbb{R}^d)$.
\end{definition}

The following theorem is the main result of this paper.

\begin{theorem}
\label{th:main} Assume that, for some $p>0$ and $q>1$, $u_0\in (\mathfrak{L}%
)_{p,q}(\mathbb{W};L_2(\mathbb{R}^d))$ and $f$, $g$ are elements of the
space $(\mathfrak{L})_{p,q}(\mathbb{W};L_2((0,T);H^{-1}(\mathbb{R}^d)))$.
Then there exist $r,\ell<0$ such that equation \eqref{eq:linG} has a unique
solution $u\in (\mathfrak{L})_{r,\ell}(\mathbb{W};L_2((0,T);H^1(\mathbb{R}%
^d)))$ and
\begin{equation}  \label{eq:main1}
\int_0^T\|u(t)\|_{r,\ell;1}^2dt \leq C\cdot \left( \|v\|_{p,q;0}^2 +\int_0^T %
\Big(\|f(t)\|_{p,q;-1}^2 +\|g(t)\|_{p,q;-1}^2\Big)\,dt \right).
\end{equation}
The number $C>0$ depends only on $\delta,C_0,p,q,r,\ell,$ and $T$.
\end{theorem}

\textbf{Proof.} The proof consists of two steps: first, we prove the result
for deterministic functions $v,f,g$ and then use linearity to extend the
result to the general case.

\emph{Step 1.} Assume that the functions $v\in L_2(\mathbb{R}^d)$, $f,g\in
L_2((0,T); H^{-1}(\mathbb{R}^d))$ are deterministic. Then $v_h=v,\, f_h=f,\,
g_h=g$, and classical theory of parabolic equations shows that, for
sufficiently small $h$, equation \eqref{eq:Str-sol} has a unique solution $%
u_h$ and the dependence of $u_h$ on $h$ is analytic.

As in the previous section, we write
\begin{equation}
u(t,x)=\sum_{\alpha\in \mathcal{J}} u_{\alpha}(t,x)\xi_{\alpha}
\end{equation}
where the coefficients $u_{\alpha}$ satisfy
\begin{equation}
\begin{split}  \label{eq:S-syst-G}
u_{(0)}(t,x)&=v(x)+\int_0^t(\mathcal{A} u_{(0)}(s,x)+f(s,x))ds, \\
u_{\epsilon_k}(t,x)&= \int_0^t\mathcal{A} u_{\epsilon_k}(s,x)ds+ \int_0^t (%
\mathcal{B} u_{(0)}(s,x)+g(s,x))m_k(s)ds, \\
u_{\alpha}(s,x)&= \int_0^t \mathcal{A} u_{\alpha}(s,x)ds+ \sum_{k} \sqrt{%
\alpha_k} \int_0^t \mathcal{B} u_{\alpha-\epsilon_k}(s,x) m_k(s)ds,\
|\alpha|>1.
\end{split}%
\end{equation}
Denote by $\Phi=\Phi_{s,t},\ t\geq s\geq 0$ the semigroup generated by the
operator $\mathcal{A}$. It follows by induction on $|\alpha|$ that
\begin{equation}
\begin{split}
u_{(0)}(t,x)& = \Phi_{t,0}v(x) +\int_0^t\Phi_{t,s}f(s)ds, \\
u_{\epsilon_k}(t,x)&= \int_0^t \Phi_{t,s}(\mathcal{B}
u_{(0)}(s,x)+g(s,x))m_k(s)ds, \\
u_{\alpha}(t,x)&=\frac{1}{\sqrt{\alpha!}} \int_0^t\int_0^{s_n} \ldots
\int_0^{s_2} \Phi_{t,s_{n}}\mathcal{B}\Phi_{s_{n},s_{n-1}} \cdots \mathcal{B}%
\Phi_{s_2,s_1}(\mathcal{B} u_{(0)}(s_1,x)+g(s_1,x)) \\
&E_{\alpha}(s_1,\ldots,s_n)ds_1\ldots ds_n, \ |\alpha|=n>1.
\end{split}%
\end{equation}
Therefore, using the usual parabolic estimates,
\begin{equation}
\int_0^T\|u_{\alpha}(t)\|_{H^1(\mathbb{R}^d)}^2dt
 \leq \frac{C^n n!}{\alpha!}
\left(\|v\|_{L_2(\mathbb{R}^d)}^2+ \int_0^T\Big( \|f(t)\|_{H^{-1}(\mathbb{R}%
^d)}^2+ \|g(t)\|_{H^{-1}(\mathbb{R}^d)}^2\Big)dt \right),
\end{equation}
and then \eqref{eq:main1} follows from \eqref{eq:pqnorm-ser}.

\emph{Step 2.} As in Step 1, existence and uniqueness of solution follows
from unique solvability of the parabolic equation \eqref{eq:Str-sol}, and it
remains to establish (\ref{eq:main1}).

Denote by $u(t,x;V,F,G,\gamma)$, $\gamma\in \mathcal{J}$, the solution of %
\eqref{eq:linG} with $v=V\xi_{\gamma}$, $f=F\xi_{\gamma}$, $g=G\xi_{\gamma}$%
. If $v=\sum_{\alpha\in \mathcal{J}}v_{\alpha} \xi_{\alpha}$, etc., then
\begin{equation}  \label{eq:sol-sum}
u(t,x)=\sum_{\gamma\in \mathcal{J}}
u(t,x;v_{\gamma},f_{\gamma},g_{\gamma},\gamma).
\end{equation}
It follows from (\ref{eq:S-syst-G}) that $u_{\alpha}(t,x;V,F,G,\gamma)=0$ if
$|\alpha|<|\gamma|$ and
\begin{equation}
\frac{u_{\alpha+\gamma}(t,x;V,F,G,\gamma)} {\sqrt{(\alpha+\gamma)!}}= \frac{%
u_{\alpha}\left(t,x;\frac{V}{\sqrt{\gamma!}}, \frac{F}{\sqrt{\gamma!}},
\frac{G}{\sqrt{\gamma!}},(0)\right)}{\sqrt{\alpha!}}.
\end{equation}
Using the results of Step 1,
\begin{equation}
\begin{split}
\int_0^T&\|u(t,\cdot;v_{\gamma}, f_{\gamma},g_{\gamma})\|_{r,\ell;1}^2dt \\
&\leq \frac{C}{\gamma!} \left(\|v_{\gamma}\|_{L_2(\mathbb{R}^d)}^2+ \int_0^T%
\Big( \|f_{\gamma}(t)\|_{H^{-1}(\mathbb{R}^d)}^2+ \|g_{\gamma}(t)\|_{H^{-1}(%
\mathbb{R}^d)}^2\Big)dt \right).
\end{split}%
\end{equation}
Now \eqref{eq:main1} follows from \eqref{eq:sol-sum} by the triangle
inequality. \endproof

%\bibliographystyle{plain}
%\bibliography{WCERef}

\begin{thebibliography}{1}

\bibitem{CM}
R.~H. Cameron and W.~T. Martin.
\newblock The orthogonal development of nonlinear functionals in a series of
  {Fourier-Hermite} functions.
\newblock {\em Ann. Math.}, 48(2):385--392, 1947.

\bibitem{HOUZ}
H.~Holden, B.~{\O}ksendal, J.~Ub{\o}e, and T.~Zhang.
\newblock {\em Stochastic Partial Differential Equations}.
\newblock Birkh\"{a}user, Boston, 1996.

\bibitem{Ito}
K.~Ito.
\newblock Multiple {Wiener} integral.
\newblock {\em J.~Math.~Soc.~Japan}, 3:157--169, 1951.

\bibitem{LR_shir}
S.~V. Lototsky and B.~L. Rozovskii.
\newblock Stochastic differential equations: a {Wiener} chaos approach.
\newblock In Yu. Kabanov, R.~Liptser, and J.~Stoyanov, editors, {\em From
  stochastic calculus to mathematical finance: the {Shiryaev} festschrift},
  pages 433--507. Springer, 2006.

\bibitem{LR_AP}
S.~V. Lototsky and B.~L. Rozovskii.
\newblock Wiener chaos solutions of linear stochastic evolution equations.
\newblock {\em Ann. Probab.}, 34(2):638--662, 2006.

\bibitem{MR_WCSIAM}
R.~Mikulevicius and B.~L. Rozovskii.
\newblock Linear parabolic stochastic PDE's and Wiener Chaos.
\newblock {\em SIAM J.~Math.~Anal.}, 292:452--480, 1998.

\bibitem{Nualart}
D.~Nualart.
\newblock {\em Malliavin Calculus and Related Topics, 2nd Edition}.
\newblock Springer, New York, 2006.

\end{thebibliography}

\end{document}